\documentclass[11pt, a4paper, twoside]{article}
\pdfoutput=1
\usepackage{lspaper} 
\newgeometry{margin=.9in}
\usepackage{tikzsymbols}
\usepackage{booktabs}
\usepackage{multirow, multicol}
\usepackage{placeins}
\usepackage[normalem]{ulem}
\usepackage{amssymb,amsmath}
\usepackage{pdflscape,afterpage}
\usetikzlibrary{matrix,arrows,decorations.pathmorphing,decorations.pathreplacing}
\usetikzlibrary{positioning,intersections,shapes.geometric,calc}
\usetikzlibrary{decorations.markings}
\RequirePackage{pxfonts}
\usepackage[toc,page]{appendix}
\usepackage{scalefnt,caption}
\usepackage[all, knot]{xy}

\renewcommand\sss[1][n]{\mathfrak{S}_{#1}}
\renewcommand\spe[1]{\rspe{#1}}

\renewcommand\D[1]{\rD{#1}}

\renewcommand\phi\varphi

\renewcommand{\theenumi}{\roman{enumi}}

\begin{document}
\setcounter{MaxMatrixCols}{19}

\title{The minimal counterexample to James's conjecture}

\author{Liron Speyer
\\\normalsize Okinawa Institute of Science and Technology\\\normalsize Onna-son, Okinawa, Japan 904-0495 \\\texttt{\normalsize liron.speyer@oist.jp}
}

\renewcommand\auth{Liron Speyer}

\runninghead{The minimal counterexample to James's conjecture}


\toptitle

\begin{abstract}
In 2017, Geordie Williamson proved the existence of counterexamples to James’s conjecture on the decomposition matrices of symmetric groups and their Hecke algebras.
The smallest counterexample detectable by Williamson's method occurs in the symmetric group $\sss$ for $n=1 \thinspace 744 \thinspace 860$, in characteristic $p=2237$.
Those detected by Williamson remain the only known counterexamples to James’s conjecture.
In this work, we calculate an explicit new counterexample, occurring in the principal block of the Hecke algebra $\hhh[24]$ when $q$ is a primitive fourth root of unity, and give explicit graded decomposition numbers in this case.
This is the minimal rank counterexample for $e\neq 2$.
\end{abstract}

\section{Introduction}

One of the most famous problems in modular representation theory is that of determining the dimensions of simple modules of the symmetric groups $\sss$ and their associated Hecke algebras $\hhh$.
This problem is mostly commonly attacked by considering the \emph{decomposition numbers} (the composition factor multiplicities of the classical (and quantised) Specht modules) which are structurally more accessible -- for example they are preserved by Morita equivalence.
The famous LLT conjecture of Lascoux--Leclerc--Thibon~\cite{LLT} -- later verified by Ariki~\cite{ariki96} -- provided a complete and beautiful algorithmic answer to this question for the Hecke algebra $\hhh$ over the complex field. 
Over fields of finite characteristic, the enduring fascination of this problem is exemplified by the recent conjecture of Lusztig--Williamson, as presented in Williamson's ICM plenary address~\cite{LusztigWilliamsonBilliards}.

Inspired by the theory of finite groups, Gordon James conjectured in~\cite{j90} that the decomposition matrices (and hence dimensions of simple modules) of the Hecke algebras should be independent of the characteristic of the field, in the so-called `abelian defect' setting.
To the astonishment of everyone in the field, Geordie Williamson disproved James's conjecture in \cite{w17}.
Williamson proved that decomposition numbers encode the arithmetic complexity of Fibonacci primes; 
this discovery effectively buried the hope that group-theoretic restrictions like abelian defect could guarantee a predictable, combinatorial answer to the problem of simple module dimensions.
We refer to \cite{chrisbook} for an excellent and broadly accessible treatment of Williamson's results, entirely in the language of diagrammatic algebra.

Williamson proved the existence of a counterexample to James's conjecture in the principal block of $\bbf \sss$ for $n = 1 \thinspace 744 \thinspace 860$, and $p = 2237$, which was the smallest known counterexample to date before the current paper.
There, Williamson writes 
\begin{quote}
\em The size of this number is a relic of our method (in particular the fact that we cannot say anything about $p$-restricted weights).
It is an important question as to where the first counter-examples occur.
\end{quote}
Here, we answer this question, giving the minimal counterexample to James's conjecture.
It occurs in the principal block of $\hhh[24]$ when $e=4$ and $p=7$; this is a weight 6 block.
Moreover, whilst Williamson's methods only show the \emph{existence} of counterexamples, we are able to explicitly calculate the exact graded decomposition numbers which nullify James's conjecture.

\begin{thm}
Let $e=4$, $\la=(5,4,3^2,2,1^7)$, and $\mu=(1^{24})$.
Then in characteristic $0$, the graded decomposition number $d_{\la,\mu} (v)$ is equal to $3v^2$, while in characteristic $7$ it is equal to $3v^2+1$.
\end{thm}

The occurrence of a counterexample in such small rank has surprised experts that we have shared our results with prior to the completion of this paper.

\begin{ack}
The author is partially supported by JSPS Kakenhi grant number 23K03043.
This material is based upon work supported by the National Science Foundation under Grant No. DMS-1929284 while the author was in residence at the Institute for Computational and Experimental Research in Mathematics in Providence, RI, during the Categorification and Computation in Algebraic Combinatorics semester program.
We thank Chris Bowman for providing sustained motivation for this project during the course of the program, Robert Muth for helpful discussions, and Jon Brundan for a helpful suggestion.
We are grateful for the help and support provided by the Scientific Computing and Data Analysis section of Core Facilities at OIST.
We thank Matt Fayers, whose GAP code for KLR algebras and Specht modules forms the basis of some of our computations, especially for computing entries of Gram matrices, and we thank Andrew Mathas for independently verifying the counterexamples we have found.
Ultimately, the part of our GAP code that handles the computations of Gram matrices in order to determine decomposition numbers comes from our work in type $C$ KLR algebras, joint with Chung and Mathas~\cite{cms25}. Finding our counterexample to James's conjecture would not have been possible without this, and we thank Andrew Mathas for explaining how to perform the necessary Gram matrix computations, which is also a key part of our GAP code for computing decomposition numbers.
\end{ack}

\section{Background}\label{sec:background}

Throughout, we let $\bbf$ denote an algebraically closed field of characteristic $p\geq 0$.
All our modules are left modules.

\subsection{Hecke algebras and their decomposition matrices}\label{subsec:Hecke}

We denote by $\hhh$ the type $A$ Iwahori--Hecke algebra, which is the associative $\bbf$-algebra with a fixed parameter $q\in \bbf^\times$, and generators $T_1,\dots,T_{n-1}$ subject to the relations
\[
(T_i -q) (T_i + 1) = 0, \quad
T_i T_{i+1} T_i = T_{i+1} T_i T_{i+1}, \quad
T_i T_j = T_j T_i \, \text{ if } \, |i-j|>1
\]
for all admissible $i,j$.
We let $e$ denote the minimal positive integer such that $1 + q + \dots + q^{e-1} = 0$, setting $e=\infty$ if no such integer exists.
An excellent introduction to the representation theory of $\hhh$ can be found in \cite{mathas}.

These algebras have Specht modules $\spe\la$ indexed by partitions $\la$ of $n$, and simple modules $\D\mu$ appearing as the simple heads of $\spe\mu$ for each $e$-restricted partition $\mu$ of $n$ -- that is for each $\mu$ such that $\mu_i - \mu_{i+1} < e$ for all $i$.
The reader should note that -- following op.~cit.~-- we work with the dual setup to that of Dipper and James's Specht modules.

The blocks of $\hhh$ can be indexed by a pair $(\rho, w)$, where $\rho$ is an $e$-core (and a partition of $n - ew$), and $w$ is the weight, or defect, of the block.
Then $\spe\la$ (respectively $\D\mu$) belongs to the block indexed by $(\rho, w)$ if and only if $\la$ (respectively $\mu$) has $e$-core $\rho$ and weight $w$.

Over $\bbc$, the decomposition numbers -- that is, the composition multiplicities $d_{\la,\mu}^p = [\spe\la : \D\mu]$ -- may be computed by the LLT algorithm~\cite{LLT}, thanks to Ariki's categorification theorem \cite{ariki96}.

\subsection{Cyclotomic KLR algebras and their graded decomposition matrices}\label{subsec:KLR}

In fact, with the advent of Khovanov--Lauda--Rouquier (KLR) algebras, we can now say more about the decomposition numbers above.
Brundan and Kleshchev~\cite{bk09} have shown that blocks of $\hhh$ are isomorphic to certain cyclotomic KLR algebras $\scrr_\beta^{\La_0}$ in type $A^{(1)}_{e-1}$ or $A_\infty$ if $e$ is finite or not, respectively.
These are $\bbz$-graded algebras, indexed by elements $\beta$ in the positive part of the corresponding root lattice, with generating set
\[
\{e(\bfi) \mid \bfi \in \{0,1,\dots,e-1\}^n\} \cup \{y_1,y_2,\dots,y_n\} \cup \{\psi_1,\psi_2, \dots,\psi_{n-1}\}
\]
subject to a long list of relations -- see \cite{kl09,rouq}.
Their simple modules and Specht modules have graded lifts, as described in \cite{bk09,bkw11,kmr}, which allow us to study the further refined \emph{graded} decomposition numbers
\[
d_{\la,\mu}^p(v) = [\spe\la : \D\mu]_v = \sum_{k\in \bbz} [\spe\la : \D\mu\langle k \rangle] \, v^k,
\]
where $M\langle k \rangle$ denotes the graded shift of the graded module $M$, so that $M\langle k \rangle_{d} = M_{d-k}$.
Note that $d_{\la,\mu}^p(1) = d_{\la,\mu}^p$.

We let $D_0 = (d_{\la,\mu}^0)_{\la,\mu}$, $D_p = (d_{\la,\mu}^0)_{\la,\mu}$ denote the ungraded decomposition matrices in characteristics $0$ and $p$, respectively, and $D_0(v) = (d_{\la,\mu}^0(v))_{\la,\mu}$, $D_p(v) = (d_{\la,\mu}^0(v))_{\la,\mu}$ denote the corresponding graded decomposition matrices.
When the partitions are ordered appropriately, these are block diagonal matrices, with blocks of the matrix naturally corresponding to blocks of $\hhh$.
It is known that there exists a lower unitriangular square matrix $A$ whose entries are non-negative integers such that
$D_p = D_0 A$,
and its graded analogue $A_v$, a lower unitriangular square matrix whose entries are Laurent polynomials in $\bbz[v + v^{-1}]$ such that
$D_p(v) = D_0(v) A_v$.
Since they are block diagonal, we will use the same notation for the block submatrices of $D_0$, $D_p$, $A$, $D_0(v)$, $D_p(v)$, and $A_v$ indexed by partitions in a fixed block of $\hhh$ when convenient.

James's original conjecture is the following, which appears in \cite[Section~4]{j90}.

\begin{conj}\label{conj:James}
Suppose $n < ep$.
Then $A$ is the identity matrix.
\end{conj}

This was slightly strengthened to a blockwise version, stated as follows.

\begin{conj}\label{conj:blockwiseJames}
For a block of $\hhh$ of weight $w < p$, the adjustment matrix $A$ is the identity matrix.
\end{conj}

One may equivalently replace $A$ with $A_v$ in the above conjectures.

By Ariki's categorification theorem~\cite{ariki96} and its graded lift due to Brundan and \linebreak Kleshchev~\cite[Corollary~5.15]{bk09}, the graded decomposition number $d_{\la,\mu}^0(v)$ is known to be equal to the coefficient of $\la$ in the canonical basis element indexed by $\mu$ in the irreducible highest weight $U_q(\widehat{\fks\fkl_e})$-module $V(\La_0)$.
In particular, there is a recursive algorithm -- the LLT algorithm~\cite{LLT} -- that computes these Laurent polynomials quite quickly.
The validity of James's conjecture would thus imply that this algorithm applies equally well in finite characteristic larger than $n/e$ (or larger than the weight $w$ of a block, in the blockwise version).

For a long time, this conjecture was suspected to be true, and indeed is known to be true for blocks of weight at most $4$, by \cite{richardswt2,faywt2,fay08wt3, fay07wt4}, as well as for the principal blocks of Hecke algebras $\hhh[5e]$ -- with the exception of two entries of the adjustment matrix when $e=4$ -- by \cite{lowwt5}.
It is also known to be true for RoCK blocks of any weight, by \cite[Corollary~3.15]{jlm}.

However, the conjecture is known to be false in general.
As stated in our introduction, Williamson~\cite{w17} has produced a counterexample in the symmetric group $\sss$ (i.e.~the case $e=p$), where $n= 1 \thinspace 744 \thinspace 860$, and $p=2237$.
Here, we answer the important question of where the first counterexample occurs.

In order to explicitly compute graded decomposition matrices $D_p(v)$ in various characteristics, we have written GAP~\cite{GAP4} code which carries out computations with formal characters of Specht modules and Gram matrix computations within certain weight spaces where necessary.
This code is explained in more detail (in type $C$) in \cite{cms25}, where Chung, Mathas, and the author first carried out such computations.

For the purposes of our counterexample, it suffices to pick out a single graded decomposition number $d_{\la,\mu}^p(v) \neq d_{\la,\mu}^0(v)$.
Although our GAP code has computed the entire graded decomposition matrix for the block where we find our counterexample, and we provide the relevant decomposition number at the end of \cref{sec:mincounterex}, it suffices to explicitly compute the corresponding Gram matrix over $\bbz$, showing that at least one of its elementary divisors is divisible by a prime $p$ larger than the weight.
Our minimal example is a counterexample for both of \cref{conj:James,conj:blockwiseJames}.
In order to describe it, we introduce the Specht modules and their Gram matrices in the next section.

\subsection{Specht modules over KLR algebras}\label{subsec:Spechts}

Let $\la$ be a partition of $n$.
For each standard $\la$-tableau $\ttt \in \std\la$ we fix a preferred choice of reduced expression $w^\ttt = s_{i_1} \dots s_{i_r}$ so that $\ttt = w^\ttt \ttt^\la$, where $\ttt^\la$ is the (row)-initial tableau.
We define the associated element $\psi^\ttt = \psi_{i_1} \dots \psi_{i_r} \in \scrr_\beta^{\La_0}$.

As a module over $\scrr_\beta^{\La_0}$, the Specht module $\spe\la$ has a homogenous presentation that may be found in \cite{kmr}.
It is generated by a single vector $v^{\ttt^\la}$,
and has a homogenous basis
\[
\{v^\ttt = \psi^\ttt v^{\ttt^\la} \mid \ttt \in \std \la\}.
\]
It will be useful for us to note that for any $\bfi \in \{0,\dots,e-1\}^n$, $e(\bfi) v^\ttt = \delta_{\bfi,\bfi^\ttt} v^\ttt$, where $\bfi^\ttt$ is the \emph{residue sequence} of $\ttt$ -- see \cite[Section~2.3]{kmr}.
In this way, we may think of the `$\bfi$-weight space' of $\spe\la$, for any $\bfi$.
Hu--Mathas~\cite{hm10} showed that these Specht modules arise as the graded cell modules of $\scrr_\beta^{\La_0}$, by providing the algebra with a graded cellular basis
\[
\{c_{\tts \ttt} = \psi^\tts e(\bfi^\la) y^\la (\psi^\ttt)^\ast \mid \tts,\ttt\in \std\la, \, \la \text{ a partition of } n\},
\]
where $*$ is the anti-involution on $\scrr_\beta^{\La_0}$ which acts trivially on the KLR generators, $\bfi^\la$ is the residue sequence of $\ttt^\la$, and $y^\la$ is a certain product of the KLR generators $y_r$ -- see \cite[Section~4.3]{hm10}.

The Specht module $\spe\la$ is equipped with a $\bbz$-valued homogeneous degree zero bilinear form $\langle - , - \rangle$ arising from the graded cellular basis above, and determined by
\[
\langle v^\tts, v^\ttt \rangle v^{\ttt^\la} = e(\bfi^\la) y^\la (\psi^\tts)^\ast v^\ttt.
\]

\section{The minimal counterexample to James's conjecture}\label{sec:mincounterex}

As remarked earlier, James's conjecture is known to be true in weight $4$ or less. It thus follows that if $e \geq 5$, James's conjecture holds for $n<25$.
Similarly, it holds for $e=4$ whenever $n\leq20$ -- with the case $n=20$ following by filling the computations missed in \cite{lowwt5} using our GAP code.
In fact, we have verified that James's conjecture holds for $e=4$ and $n\leq 23$.
The counterexample we have found in the principal block of $\hhh[24]$ is thus the smallest rank counterexample for $e\geq 4$.

By our GAP computations, we have ruled out any possible counterexamples in $\hhh$ for $e=2$ and $n \leq 14$, and for $e=3$ and $n \leq 19$.
Though we have not yet managed to rule out a counterexample for $e=2$ and $15 \leq n \leq 23$, or for $e=3$ and $20 \leq n \leq 23$, these still amount to $n > 6e$, so that our example is still of minimal rank relative to $n$.
But Mathas has further verified that, in fact, no smaller rank counterexamples exist for $e=3$, and no counterexamples exist for $e=2$ and $n \leq 17$.
For $e=2$ and $18 \leq n \leq 19$, Mathas has checked that no counterexamples exist in the $\bfi^{(1^n)}$ weight-spaces of any Specht modules, so we expect that no counterexamples exist in those ranks, either. We are grateful to Mathas for performing these checks for us, as they greatly bolster our case that we have indeed found the minimal counterexample to James’s conjecture.

Recall that we work with the `row' Specht modules of \cite{kmr}, coinciding with those of \cite{mathas}, and dual to those originally studied in the works of Dipper and James.
We now fix the 4-restricted partitions $\la = (5,4,3^2,2,1^7)$ and $\mu=(1^{24})$, and will focus on the multiplicity of $\D\mu$ in $\spe\la$, noting that these partitions are in the principal block of $\hhh[6e]$ which has weight 6.
To this end, we may check using the LLT algorithm that $d_{\la,\mu}^0(v) = 3v^2$.

We check the $\bfi^\mu$-weight space of $\spe\la$ -- i.e.~the subspace
\[
e(\bfi^\mu)\spe\la =
e(0,3,2,1,0,3,2,1,0,3,2,1,0,3,2,1,0,3,2,1,0,3,2,1) \spe\la.
\]
One may check that the graded dimension of $e(\bfi^\mu)\spe\la$ is $9v^4+43v^2+19$.
In fact, we zoom in even further, and will look at the 19-dimensional degree 0 subspace of $e(\bfi^\mu)\spe\la$.
Computing (in GAP) the corresponding $19\times19$ submatrix of the Gram matrix $(\langle v^\tts, v^\ttt \rangle)_{\tts,\ttt\in\std\la}$ yields the following matrix.
\[
\begin{bmatrix}
-2 & -1 & 1 & 0 & 0 & 0 & 0 & 0 & 0 & 0 & 0 & 0 & 0 & 0 & 0 & 0 & 0 & 0 & 0 \\
-1 & -2 & 0 & 1 & -1 & 0 & 1 & 0 & 0 & 0 & 0 & 0 & 0 & 0 & 0 & 0 & 0 & 0 & 0 \\
1 & 0 & -2 & 0 & 0 & 0 & 0 & 1 & 0 & 0 & 0 & 0 & 0 & 0 & 0 & 0 & 0 & 0 & 0 \\
0 & 1 & 0 & -2 & 1 & 0 & -2 & 0 & 0 & 0 & 0 & 0 & 0 & 0 & 0 & 0 & 0 & 0 & 1 \\
0 & -1 & 0 & 1 & -2 & 0 & 2 & 0 & 0 & 0 & 1 & 0 & -1 & 0 & 0 & 0 & 0 & 0 & 0 \\
0 & 0 & 0 & 0 & 0 & 0 & 0 & 0 & 0 & 0 & 0 & 0 & 0 & 0 & 0 & 0 & 0 & 0 & 0 \\
0 & 1 & 0 & -2 & 2 & 0 & -4 & 0 & 0 & 1 & -1 & 0 & 2 & 0 & 0 & -1 & 0 & 0 & 2 \\
0 & 0 & 1 & 0 & 0 & 0 & 0 & -2 & 0 & 0 & 0 & 0 & 0 & 1 & 0 & 0 & 0 & 0 & 0 \\
0 & 0 & 0 & 0 & 0 & 0 & 0 & 0 & 0 & 0 & 0 & 0 & 0 & 0 & 0 & 0 & 0 & 0 & 0 \\
0 & 0 & 0 & 0 & 0 & 0 & 1 & 0 & 0 & -2 & 0 & 0 & -1 & 0 & 0 & 2 & 0 & 0 & -1 \\
0 & 0 & 0 & 0 & 1 & 0 & -1 & 0 & 0 & 0 & -2 & 0 & 2 & 0 & 0 & 0 & 0 & 0 & 0 \\
0 & 0 & 0 & 0 & 0 & 0 & 0 & 0 & 0 & 0 & 0 & 0 & 0 & 0 & 0 & 0 & 0 & 0 & 0 \\
0 & 0 & 0 & 0 & -1 & 0 & 2 & 0 & 0 & -1 & 2 & 0 & -4 & 0 & 0 & 2 & 0 & 0 & -1 \\
0 & 0 & 0 & 0 & 0 & 0 & 0 & 1 & 0 & 0 & 0 & 0 & 0 & -2 & 0 & 0 & 0 & 0 & 0 \\
0 & 0 & 0 & 0 & 0 & 0 & 0 & 0 & 0 & 0 & 0 & 0 & 0 & 0 & 0 & 0 & 0 & 0 & 0 \\
0 & 0 & 0 & 0 & 0 & 0 & -1 & 0 & 0 & 2 & 0 & 0 & 2 & 0 & 0 & -4 & 0 & 0 & 2 \\
0 & 0 & 0 & 0 & 0 & 0 & 0 & 0 & 0 & 0 & 0 & 0 & 0 & 0 & 0 & 0 & 0 & 0 & 0 \\
0 & 0 & 0 & 0 & 0 & 0 & 0 & 0 & 0 & 0 & 0 & 0 & 0 & 0 & 0 & 0 & 0 & 0 & 0 \\
0 & 0 & 0 & 1 & 0 & 0 & 2 & 0 & 0 & -1 & 0 & 0 & -1 & 0 & 0 & 2 & 0 & 0 & -4
\end{bmatrix}
\]
The elementary divisors of this matrix are $\alpha_1 = \alpha_2 = \dots = \alpha_{12} = 1$, $\alpha_{13} = 14$, and $\alpha_{14} = \dots = \alpha_{19} = 0$.
In other words, the rank of the matrix decreases by one in characteristic $7$ (and also in characteristic $2$).
In particular, this means that the dimension of the degree $0$ subspace of $e(\bfi^\mu)\D\la$ decreases by one in characteristic $7$, which violates \cref{conj:James,conj:blockwiseJames}.\footnote{The GAP code that constructs this weight space, generates its degree-zero tableaux and the corresponding reduced expressions, and reproduces the $19\times19$ Gram matrix and its invariant factors is supplied as ancillary material with the arXiv version of this paper.
Instructions are given in \texttt{README.md}; running \texttt{james\_counterexample.g} performs the calculation.}

In fact, we have calculated in GAP that $d_{\la,\nu}^7(v) = d_{\la,\nu}^0(v)$ for all $\nu \neq \mu$, but that $d_{\la,\mu}^7(v) = 3v^2 +1 \neq d_{\la,\mu}^0(v)$.
Moreover, this gives us that the entry in the $(\la,\mu)$-position of the adjustment matrix is a $1$, and therefore there are many other $\gamma$ for which $d_{\gamma,\mu}^7(v) \neq d_{\gamma,\mu}^0(v)$ 

\begin{rmk}
\begin{enumerate}
\item
The above block has one further non-trivial entry in the adjustment matrix in characteristic $7$ -- namely another $1$ in the $((7,6,5,3,2,1),(3^8))$-position.
Here, $\D{(3^8)}$ is the trivial module.
In fact, this entry can be determined from the one we have explained above, using the fact that in the adjustment matrix, $A_{\la\mu} = A_{\la^\diamond \mu^\diamond}$, by \cite[Lemma~4.2]{fay08wt3}, where $\diamond$ is the Mullineux involution on $e$-restricted partitions of $n$.

\item
We are grateful to Andrew Mathas for independently verifying our computation of the Gram matrix for the weight space we have computed above.
After verifying our result, Mathas has implemented a search for more counterexamples by looking at the $\bfi^{(1^n)}$-weight spaces of Specht modules and computing the Gram matrices, and we have implemented a similar search in GAP.
Mathas's code -- computing the products necessary to determine Gram matrix entries by working with seminormal forms -- appears to be much quicker than our own, and he has now found a wealth of these examples.
Namely, there are many more for larger ranks in weight 6 blocks and characteristic 7.
Perhaps the most important of these, which we have also managed to verify, is one in the principal block of $\sss[44]$ in characteristic $7$.
Though there are no smaller rank counterexamples in symmetric groups that occur in these $\bfi^{(1^n)}$-weight spaces of Specht modules, we have not ruled out counterexamples in ranks $n$ for $36 \leq n \leq 43$ when $p=7$.
Nonetheless, we expect that this $\sss[44]$ example is the minimal one in the symmetric group case, and that these $\bfi^{(1^n)}$-weight spaces of Specht modules detect the smallest rank counterexamples to James's conjecture.
Another important counterexample Mathas finds, and that we have verified -- also in a weight 6 block in characteristic 7 -- occurs in $\hhh[37]$ when $e=6$.
Strictly speaking, since \cite{j90} is concerned with the finite group $GL_n(q)$, it assumes throughout that the Hecke parameter $q$ is a prime power, and in particular an integer, and therefore an element of the prime subfield $\bbf_p$.
While this condition is not usually considered in the literature on James's conjecture, this $\hhh[37]$ counterexample appears to be the first for which such a $q$ exists, therefore disproving even the strictest interpretation of James's conjecture.
We thank Meinolf Geck for pointing out this extra condition to us, and asking about the existence of a counterexample (with $q\neq 1$) in such cases.

It is worth noting that Mathas has also found counterexamples in higher weights and characteristics.
His examples will appear in \cite{MathasJamesCounters}.
\end{enumerate}
\end{rmk}

We end with a conjecture, based on the computations made in search for our counterexample.

\begin{conj}
James's conjecture holds for all weight $5$ blocks of $\hhh$.
\end{conj}

Given our lack of counterexamples in quantum characteristics 2 and 3, it would also be interesting to determine the veracity of James’s conjecture in these cases.

\bibliographystyle{lspaper} 
\addcontentsline{toc}{section}{\refname}
\bibliography{master}

\end{document}